\documentclass[10pt]{amsart}
\usepackage{amscd,amssymb}
\usepackage[dvips]{graphics}

\newtheorem{thm}{Theorem}[section]
\newtheorem{lemma}[thm]{Lemma}

\newtheorem{defn}[thm]{Definition}

\newcommand{\R}{\mathbb R}

\renewcommand{\o}{\operatorname}

\newcommand{\oto}{\o{O}(2,1)}
\newcommand{\soto}{{\o{SO}(2,1)^0}}

\newcommand{\rto}{\R^{2,1}}

\newcommand{\SLR}{\o{SL}(2,{\mathbb R})}
\newcommand{\PSLR}{\o{PSL}(2,{\mathbb R})}
\newcommand{\slr}{\mathfrak{sl}(2,\R)}

\newcommand{\E}{\mathbb M^{2,1}} 
\newcommand{\Gg}{{\mathsf{G}}} 
\newcommand{\Go}{{\Gg}} 
\newcommand{\newalpha}{\tilde{\alpha}}
\newcommand{\iso}{\o{Aff}(\E)}

\newcommand{\hyp}{{\mathbb H}^2}
\newcommand{\psl}{\o{Isom}(\hyp)} 

\newcommand{\LL}{\mathbb L}

\newcommand{\vx}{{\mathsf x}}
\newcommand{\vy}{{\mathsf y}}

\newcommand{\vv}{{\mathsf v}}

\newcommand{\lv}{V}
\newcommand{\lw}{W}
\newcommand{\lbasis}{E}

\newcommand{\lx}{X}
\newcommand{\lf}{F}

\renewcommand{\lg}{\mathfrak{g}}

\newcommand{\ldot}[1]{\mathbb B(#1)}

\newcommand{\ga}{\alpha}
\renewcommand{\gg}{\gamma}

\newcommand{\G}{\Gamma_0}

\newcommand{\hp}{\mathbb{H}^2}

\newcommand{\ZZ}{Z^1(\G,\R^{2,1})}
\newcommand{\HH}{H^1(\G,\R^{2,1})}
\newcommand{\ZZg}{Z^1(\G,\lg)}

\newcommand{\surf}{\Sigma}
\newcommand{\Hom}{\mathsf{Hom}}

\begin{document}

\title{Stretching Three-Holed Spheres and the Margulis Invariant}

\author[Charette]{Virginie Charette}
    \address{D\'epartement de math\'ematiques\\ Universit\'e de Sherbrooke\\
    Sherbrooke, Quebec, Canada}
    \email{v.charette@usherbrooke.ca}

\author[Drumm]{Todd A. Drumm}
    \address{Department of Mathematics\\ Howard University \\
    Washington, DC }
    \email{tdrumm@howard.edu}

\author[Goldman]{William Goldman}
    \address{Department of Mathematics\\ University of Maryland\\
    College Park, MD}
    \email{wmg@math.umd.edu}

\date{\today}

\maketitle

\section{Introduction}\label{sec:intro}

This note concerns an application of the emerging theory of
{\em complete flat Lorentz  $3$-manifolds\/}
to hyperbolic geometry on surfaces.

We shall apply our forthcoming paper~\cite{CDG}
to prove the following simple result:

\begin{thm}\label{thm:lengths}

Let $\surf$ be a three-holed sphere.  Consider a one-parameter
family $\surf_t$ of marked hyperbolic structures on $\surf$.
For each $\gamma\in \pi_1(\Sigma)$ denote the length of the closed
geodesic corresponding to $\gamma$ by $\ell(\gamma)$.

Suppose that for each $\partial_i$ corresponding to a component of
$\partial\Sigma$,
\begin{equation*}
\frac{d\ell(\partial_i)}{dt}\bigg|_{t=0} > 0.
\end{equation*}
Then for every essential closed curve $\gamma$,
\begin{equation*}
\frac{d\ell(\gamma)}{dt}\bigg|_{t=0} > 0.
\end{equation*}
\end{thm}

A {\em complete flat Lorentz $3$-manifold\/} is a
geodesically complete Lorentzian $3$-manifold
of zero curvature. Such a manifold is a quotient
$M = \E/\Gamma$  of
$(2+1)$-dimensional Minkowski space $\E$ by a discrete
group  $\Gamma$ of isometries acting properly and freely on $\E$.
Recall that
{\em  $(2+1)$-dimensional Minkowski space\/} is a geodesically complete
simply connected Lorentzian manifold  of zero curvature.

Alternatively, $\E$ is a $3$-dimensional affine space,
together with a quadratic form of signature $(2,1)$
on the vector space of translations.
We shall call such an inner product space a {\em Lorentzian vector
space.\/} Every such Lorentzian vector space is isomorphic to
$\R^3$ with
inner product
\begin{equation*}
\ldot{\vx,\vy} \;=\; x_1y_1+x_2y_2-x_3y_3,
\end{equation*}
where
\begin{equation*}
\vx= \bmatrix x_1 \\ x_2 \\ x_3\endbmatrix , \vy= \bmatrix y_1 \\
y_2 \\ y_3\endbmatrix.
\end{equation*}
We denote this Lorentzian vector space by $\rto$.
The group of automorphisms of $\rto$ is the orthogonal group
$\oto$.

The deformation theory of hyperbolic structures on surfaces
intimately relates to discrete groups of isometries of $\E$.
Every quotient of $\E$ by a discrete group $\Gamma$ of
isometries acting properly determines a complete noncompact hyperbolic
surface $\Sigma$.
In \cite{CDG}, quotients are classified when the corresponding
hyperbolic surface $\Sigma$ is homeomorphic to a $3$-holed sphere.
Here we describe this classification in terms
of deformations of hyperbolic structures on the three-holed sphere.

The interplay between the two deformation theories owes to
the identification (first exploited in this context in \cite{GM} and
\cite{G0}) of the Lie algebra of Killing vector fields on the
hyperbolic plane $\hyp$ as a Lorentzian vector space.  Infinitesimal
deformations of a hyperbolic surface $\surf$ correspond to {\em affine
deformations\/} of the holonomy representation $\pi_1(\surf)
\xrightarrow{\rho} \psl$.  By \cite{FG}, complete flat Lorentz
$3$-manifolds $M$ fall into two distinct types.  The first type arises
when $\pi_1(M)$ is solvable, in which the classification is a simple
exercise in linear algebra.  The second type, which we call {\em
nonelementary\/}, arises from an affine deformation $\rho$ of a
homomorphism $\pi_1(M) \xrightarrow{\rho_0} \psl\cong \soto$
satisfying the following conditions:

\begin{itemize}
\item $\rho_0$ is injective;

\item The image $\Gamma_0$ of  $\rho_0$
is a discrete subgroup of $\soto$;

\item $\rho$ defines a {\em proper action\/} of $\pi_1(M)$ on $\E$.

\end{itemize}

By Mess~\cite{Mess}, the hyperbolic surface $\surf :=
\hyp/\Gamma_0$ is noncompact. Thus isometry classes of
nonelementary complete flat Lorentz $3$-manifolds $M$ correspond
to conjugacy classes of affine deformations $\rho$ of discrete
embeddings of $\pi_1(M)$ in $\psl \cong \soto$, defining {\em
proper\/} actions of $\pi_1(M)$ on $\E$. Henceforth we refer to
such affine deformations $\rho$ as {\em proper deformations.}

Determining which affine deformations are proper
is a fundamental and difficult problem.
When $\pi_1(\surf)\xrightarrow{\rho_0} \psl$ embeds $\pi_1(\surf)$
in a Schottky group, \cite{GLM} provides criteria in terms of an
extension of an invariant discovered by Margulis~\cite{Margulis1,Margulis2}.
One of the main results of \cite{CDG} is that these criteria take a particularly
simple form when $\surf$ is a $3$-holed sphere.

Nonelementary flat Lorentz 3-manifolds behave like hyperbolic
surfaces in many ways. For example, if $\gamma\in \pi_1(M)$ does
not correspond to a cusp of $\surf$, then $\gamma$ corresponds to
a closed geodesic in $\surf$, and we denote its length by
$\ell(\gamma)$. In $M$, $\gamma$ corresponds to a closed geodesic
with respect to the induced flat Lorentz metric on $M$. This
geodesic is spacelike and has a well-defined positive length.

Margulis defined, for {\em any\/} affine deformation $\rho$,
a class function
\begin{equation*}
\{ \gamma\in \pi_1(M) \mid \rho_0(\gamma)\text{ is hyperbolic} \}
\xrightarrow{\alpha_\rho} \R
\end{equation*}
with many remarkable properties.

For example, $\alpha(\gamma) \neq 0$ if and only if the cyclic group
$\langle \rho(\gamma)\rangle$ acts freely on $\E$.
When $\rho$ is proper, then $\vert\alpha(\gamma)\vert$
is the length of the closed geodesic in $M$ corresponding to $\gamma$.
Furthermore in this case, either the values of $\alpha$ are all positive
or all negative.

Affine deformations of the holonomy representation
$\pi_1(\Sigma)\xrightarrow{\rho_0}\psl$ of a complete hyperbolic
surface $\Sigma$ correspond to infinitesimal deformations of the
hyperbolic manifold $\Sigma$ as follows.  As first observed by
Weil~\cite{Weil1,Weil2} (see also Raghunathan~\cite{Raghunathan}),
infinitesimal deformations of the geometric structure on $\Sigma$
correspond to vectors in the Zariski tangent space to
\begin{equation*}
\Hom\big(\pi_1(\Sigma),\psl\big)
\end{equation*}
at $\rho_0$. Weil computed this tangent space
as the cohomology
\begin{equation*}
H^1(\Sigma,\lg)
\end{equation*}
with coefficients in the local system corresponding to the Lie
algebra $\lg\cong\slr$ of Killing vector fields, that is, {\em
infinitesimal isometries \/} of $\hyp$. This local system identifies
with the $\pi_1(M)\cong\pi_1(\Sigma)$-module corresponding to the
Lorentzian vector space $\rto$ with the action given by $\rho_0$.

For example suppose that
$u\in Z^1\big(\pi_1(\Sigma),\rto\big)$ is a cocycle, the
corresponding affine deformation is explicitly given as:
\begin{equation*}
\rho(\gamma): x \longmapsto \gamma(x) + u(\gamma).
\end{equation*}
Suppose that $\rho_t$ is a smooth path in
$\Hom\big(\pi_1(\Sigma),\psl\big)$, of holonomies of hyperbolic structures,
 with geodesic length functions
\begin{equation*}
\pi_1(\Sigma) \xrightarrow{\ell_t} \R
\end{equation*}
and velocity vector $u$.

In this interpretation, Margulis's invariant $\alpha(\gamma)$ is
the rate of change of the geodesic length function
$\ell_t(\gamma)$ under the smooth path $\rho_t$ above:
\begin{equation*}
\alpha(\gamma) =
\frac{d\ell_t(\gg)}{dt}\bigg|_{t=0}.
\end{equation*}
Thus, if $\alpha(\gamma) > 0$ for $\gamma\neq 1$, then
the length of $\gamma$ increases to first order under the
deformation of hyperbolic structures corresponding to $\rho_t$.

Let $\surf$ be a noncompact surface with a complete hyperbolic
metric. If there exists a simple closed geodesic $\gamma$ bounding
a noncompact part of $\surf$ homeomorphic to a cylinder, that part
of the surface is called an {\em end} with corresponding geodesic
$\gamma$. If there is a simple closed curve bounding a noncompact
part of $\surf$ homeomorphic to a cylinder but with finite area,
that part of the surface is called a {\em cusp}.   Ends have
hyperbolic holonomy and cusps, parabolic holonomy.

We can naturally associate a length to each end of $\surf$,
namely,  the length of its corresponding closed geodesic.  Extend
this notion of length to a cusp by declaring a cusp to be of
length zero.  So given a path of holonomy representations, we can
say that the cusp {\em lengthens} if its image along the path
deforms to an end.

In this language, we may restate Theorem~\ref{thm:lengths} as follows.

\begin{thm}\label{thm:lengths2}
Let $\surf$ be a complete surface homeomorphic to a three-holed sphere with hyperbolic
structure induced by the holonomy representation
\begin{equation*}
\rho:\pi_1(\surf)\longrightarrow\psl.
\end{equation*}
Suppose $\rho_t$  is a path of holonomy representations such that
$\rho_0=\rho$.

If the lengths of the three components $\partial_1,\partial_2,\partial_3$ of
$\partial\surf$ are
increasing along $\rho_t$, then up to first order, the length of
every closed geodesic is increasing.
\end{thm}

Using an argument inspired by Thurston~\cite{Thurston}, we will in
fact prove the following  extension of Theorem~\ref{thm:lengths2}.

\begin{thm}\label{thm:thurston}
Let $\surf$ be a complete surface homeomorphic to a three-holed
sphere with boundary components
$\partial_1,\partial_2,\partial_3$.  Let $\rho_t$ be a path of
holonomy representations.  If the lengths of the $\partial_i$ are
increased (respectively  not decreased), then every closed
geodesic on $\surf$ is increased (respectively  not decreased).
\end{thm}

\section{Background}

\subsection{The geometry of $\R^{2,1}$}

Let $\E$ denote {\em Minkowski (2+1)-space}, the
three-dimensional affine space with the following additional
structure.  Its associated vector space of directions is
\begin{equation*}
\R^{2,1}=\{ p-q\mid p,q\in\E\} .
\end{equation*}
This vector space is isomorphic to $\R^3$ as a vector space with
the standard {\em Lorentzian inner product}:
\begin{equation*}
\ldot{\vx,\vy} \;=\; x_1y_1+x_2y_2-x_3y_3,
\end{equation*}
where
\begin{equation*}
\vx= \bmatrix x_1 \\ x_2 \\ x_3\endbmatrix , \vy= \bmatrix y_1 \\
y_2 \\ y_3\endbmatrix.
\end{equation*}
A non-zero vector $\vx$ is said to be {\em null} (respectively
{\em timelike}, {\em spacelike}) if $\ldot{\vx,\vx}=0$
(respectively $\ldot{\vx,\vx}<0$, $\ldot{\vx,\vx}>0$).

\subsection{Lorentzian transformations and affine deformations}

Let $\iso$ denote the group of all affine transformations that
preserve the Lorentzian inner product on the space of directions;
$\iso$ is isomorphic to $\oto\ltimes\R^{2,1}$.  We shall restrict
our attention to those transformations whose linear parts are in
$\soto$, thus preserving orientation and time-orientation.

Denote projection onto the {\em linear part} of an affine
transformation by:
\begin{equation*}
\iso\xrightarrow{\LL}\oto .
\end{equation*}
Recall that the upper nappe sheet of the hyperboloid of
unit-timelike vectors in $\rto$ is a model for the hyperbolic
plane $\hp$.  The resulting isomorphism between $\soto$ and $\psl$
gives rise to the following terminology.  (Consult~\cite{GM}, for
example, for an explicit isomorphism.)
\begin{defn}
Let $g\in\soto$ be a nonidentity element;
\begin{itemize}
\item $g$ is {\em hyperbolic} if it has three, distinct real
eigenvalues;

\item $g$ is {\em parabolic} if its only eigenvalue is 1;

\item $g$ is {\em elliptic} otherwise.
\end{itemize}
We also call $\gamma\in\iso$ {\em hyperbolic} (respectively  {\em
parabolic}, {\em elliptic}) if its linear part $\LL(\gg)$ is
hyperbolic (respectively  parabolic, elliptic).
\end{defn}

Let $\G\subset \oto$ be a subgroup. An {\em affine deformation}
of $\G$ is a representation
\begin{equation*} \rho: \G \longrightarrow \iso.
\end{equation*}
For $g\in\G$, write
\begin{equation*}
\rho(g) (x) = g(x) + u(g)
\end{equation*}
where $u(g)\in\rto$. Then $u$ is a cocycle of $\G$ with
coefficients in the $\G$-module $\rto$ corresponding to the
linear action of $\G$. In this way affine deformations of $\G$
correspond to cocycles in $\ZZ$ and translational conjugacy
classes of affine deformations correspond to cohomology classes in
$\HH$.

By extension, if $\G=\rho_0(\pi_1(\surf))$, we will call $\rho$ an
affine  deformation of the holonomy representation $\rho_0$.

\subsection{The Lie algebra $\slr$ as $\rto$}

The Lie algebra $\slr$ is the tangent space to $\PSLR$ at the
identity and consists of the set of traceless $2\times 2$
matrices. The three-dimensional vector space has a natural inner
product, the Killing form, defined to be
\begin{equation}
\left< \lv, \lw \right> = \frac{1}{2} \o{Tr} ( \lv\cdot\lw).
\end{equation}
A basis for $\slr$ is given by
\begin{equation}
 \lbasis_1=\begin{bmatrix} 1&0 \\ 0 & -1 \end{bmatrix}, \,
\lbasis_2=\begin{bmatrix} 0&1 \\ 1 & 0 \end{bmatrix}, \,
\lbasis_3=\begin{bmatrix} 0& 1 \\ -1 & 0 \end{bmatrix}.
\end{equation}
Evidently, $\langle \lbasis_1, \lbasis_1 \rangle \,=\, \langle
\lbasis_2,\; \lbasis_2 \rangle \,=\, 1$,  $\langle \lbasis_3,
\lbasis_3 \rangle \,=\,-1 $ and $\left< \lbasis_i, \lbasis_j
\right> =0$ for $i\neq j$. That is, $\slr$ is isomorphic to $\rto$
as a vector space
\begin{equation*}
\left\{ \vv = \begin{bmatrix} x\\y\\z \end{bmatrix} \right\}
\leftrightarrow \{ x\lbasis_1 + y\lbasis_2 + z\lbasis_3 =\lv \}.
\end{equation*}
The adjoint action of $\PSLR$ on $\slr$:
\begin{equation*}
g(\lv)=  g\lv g^{-1}
\end{equation*}
corresponds to the linear action of $\soto$ on $\rto$.

Using
these identifications, set:
\begin{align*}
\Go &\cong \PSLR \cong\soto \\
\lg & \cong \slr\cong\rto.
\end{align*}

\subsection{The Margulis invariant}

The Margulis invariant is a measure of an affine
transformation's signed Lorentzian displacement in $\E$. Originally
defined by Margulis for hyperbolic
transformations~\cite{Margulis1,Margulis2}, it admits an extension to parabolic
transformations~\cite{CharetteDrumm1}.

Let $g\in\Go$ be a non-elliptic element.  Lift $g$ to a representative in $\SLR$;
then the following element of $\lg$ is a $g$-invariant vector which is independent of choice of lift:
\begin{equation*}
\lf_g=\sigma(g)\left(g-\frac{\o{Tr}(g)}{2}I\right)
\end{equation*}
where $\sigma(g)$ is the sign of the trace of the lift.

Now let $\G\subset\Go$ such that every element other than the
identity is non-elliptic. Let $\rho$ be an affine deformation of
$\G$, with corresponding $u\in\ZZ$. Given the above identification
$\lg\cong\rto$, we may also write $u\in\ZZg$. We define the {\em
non-normalized Margulis invariant} of $\rho(g)\in\rho(\G)$ to be:
\begin{equation}\label{eq:alphadef}
\newalpha_\rho(g)= \langle u(g), \lf_g \rangle.
\end{equation}

(In~\cite{CharetteDrumm1}, the non-normalized invariant is a functional on a fixed line, rather than a value.)

If $\rho(g)$ happens to be hyperbolic, then $\lf_g$ is spacelike
and we may replace it by the unit-spacelike vector:
\begin{equation*}
\lx^0_g=\frac{2\sigma(g)}{\sqrt{\o{Tr}(g)^2-4}}\left(g-\frac{\o{Tr}(g)}{2}I\right)
\end{equation*}
obtaining the {\em normalized Margulis invariant}:
\begin{equation}
\alpha_\rho(g) = \langle u(g), \lx^0_g \rangle .
\end{equation}

In Minkowski space, $\alpha_\rho(g)$ is the {\em signed Lorentzian
length} of a closed geodesic in $\E/\langle \rho(g) \rangle$~\cite{Margulis1,Margulis2}.

As a function of word length in the group $\G$, normalized
$\alpha_\rho$ behaves better than  non-normalized
$\newalpha_\rho$. Nonetheless, the sign of $\newalpha_\rho(g)$ is
well defined and is equal to that of $\alpha_\rho(g)$. So we may
extend the definition of $\alpha_\rho$ to parabolic $g$, for
instance by declaring that $\lf_g=\lx^0_g$.

\begin{thm}\label{thm:tool}\cite{CDG}
Let $\G$  be a Fuchsian group whose corresponding hyperbolic
surface $\Sigma$ is homeomorphic to a three-holed sphere. Denote
the generators of $\G$ corresponding to the three components of
$\partial\Sigma$ by $\partial_1, \partial_2, \partial_3$.  Let
$\rho$ be an affine deformation of $\G$.

If
$\alpha_\rho(\partial_i) $ is positive
(respectively, negative, nonnegative, nonpositive) for each $i$ then
for all $\gamma\in\G\setminus\{ 1 \}$, $\alpha_\rho(\gamma)$ is positive (respectively, negative,
nonnegative, nonpositive).
\end{thm}

The proof of Theorem~\ref{thm:tool} relies upon showing that the
affine deformation $\rho$ of the Fuchsian group $\G$ acts properly
on $\E$. By a fundamental lemma due to Margulis
\cite{Margulis1,Margulis2} and extended in \cite{CharetteDrumm1},
if $\rho$ is proper, then $\alpha_\rho$
applied to every element has the same sign. Moreover,
 \begin{itemize}
 \item if $\alpha_\rho(\partial_1) =0$ and
$\alpha_\rho(\partial_2), \alpha_\rho(\partial_3)  >0 $ then
specifically  $\alpha_\rho(\gamma)=0$ only if $\gamma\in \langle
\partial_1\rangle$, and
 \item if $\alpha_\rho(\partial_1) =\alpha_\rho(\partial_2)=0$ and
 $\alpha_\rho(\partial_3)>0 $ then specifically $\alpha_\rho(\gamma)=0$
 only if $\gamma\in \langle \delta_1, \delta_2 \rangle$.
\end{itemize}

\section{Length changes in deformations}

As we pointed out in the Introduction, an {\em affine} deformation
of a  holonomy representation corresponds to an {\em infinitesimal
deformation} of the holonomy representation, or a tangent vector
to the holonomy representation.  In this section, we will further
explore this correspondence, relating the affine Margulis
invariant to the derivative of length along a path of holonomy
representations.  We will then prove Theorems~\ref{thm:lengths}
and~\ref{thm:thurston} by applying Theorem~\ref{thm:tool}, which
characterizes proper deformations in terms of the Margulis
invariant, to the study of length changes along a path of holonomy
representations.  We will close the section with some explicit
computations of first order length changes.

Let $\rho_0: \pi_1(\surf)\rightarrow \G\subset \Go$ be a holonomy representation 
and let $\rho:\G\rightarrow\iso$ be an affine deformation of $\rho_0$, with corresponding cocycle $u\in\ZZg$.

The affine deformation $\rho$ induces a path of holonomy representations $\rho_t$ as follows:
\begin{align*}
\rho_t:  \pi_1(\surf) & \longrightarrow \Go \\
 \gamma & \longmapsto \exp(tu(g))g,
 \end{align*}
where $g=\rho_0(\gamma)$, and $u$ is the tangent vector to this
path at $t=0$.   Conversely, for any path of representations
$\rho_t$
\begin{equation*}
\rho_t(\gamma)=\exp(tu(g)+O(t^2))g,
\end{equation*}
where $u\in\ZZg$ and $g=\rho_0(\gamma)$.

Suppose $g$ is hyperbolic.  Then the length of the corresponding closed geodesic in $\surf$ is
 \begin{equation*}
 l(g)=2\cosh^{-1}\left(\frac{|\o{Tr}(\tilde{g})|}{2}\right),
 \end{equation*}
where $\tilde{g}$ is a lift of $g$ to $\SLR$.
With $\rho,~\rho_t$ as above and $\rho_0(\gamma)=g$, set:
 \begin{equation*}
 l_t(\gamma)=l(\rho_t(\gamma)).
 \end{equation*}

Since the Margulis invariant of $\rho$ can also be seen to  be a
function of its corresponding cocycle $u$, we may write:
\begin{equation*}
\alpha_u(g):=\alpha_\rho(g).
\end{equation*}
Consequently:
 \begin{equation*}
 \frac{d}{dt}\Big|_{t=0}l_t(\gamma)=\frac{\alpha_u(g)}{2},
 \end{equation*}
so we may interpret $\alpha_u$ as the change in length of an affine deformation, up to first order~\cite{GM,G0}.

Although $l_t(\gamma)$ is not differentiable at 0 for parabolic
$g$,

  \begin{equation*}
  \frac{d}{dt}\Big|_{t=0}\frac{\sigma(g)}{2}\o{Tr}(\rho_t(\gamma))=\newalpha_u(g).
  \end{equation*}
 Thus Theorem~\ref{thm:lengths} simply
reinterprets Theorem~\ref{thm:tool}.

\begin{proof}[Proof of Theorem~\ref{thm:thurston}]

Let $\rho_t$, $-\epsilon\leq t\leq \epsilon$ be a path of holonomy
representations.  Since we assume the  boundary components to be
lengthening, they must have hyperbolic holonomy on
$(-\epsilon,\epsilon)$.

Suppose there exists $\gamma\in\pi_1(\surf)$ and
$T\in(-\epsilon,\epsilon)$ such that the  length of
$\rho_t(\gamma)$ decreases in a neighborhood of $T$.  Let
$u=u_T\in\ZZg$ be a cocycle tangent to the path at $T$; then
\begin{equation*}
\alpha_{u_T}(\gamma)<0.
\end{equation*}
Theorem~\ref{thm:tool} implies that for some $i=1,2,3$:
\begin{equation*}
 \alpha_{u_T}(\partial_i))<0.
 \end{equation*}
 but then the length of the corresponding end must decrease,
 contradicting the hypothesis.

\end{proof}

\subsection{Deformed hyperbolic transformations}

In this and the next paragraph, we explicitly compute the trace of
some deformations, to understand first order length changes.

Let $g\in\SLR$ be a hyperbolic element, thus a lift of a
hyperbolic  isometry of $\hp$.  Given a tangent vector in
$\lv\in\slr$, consider the following two actions on $\SLR$:

\begin{equation}\label{eq:lieaction}
\pi_{\lv}:g\rightarrow \o{exp}(\lv)\cdot g,
\end{equation}
and
\begin{equation}\label{eq:altlieaction}
\pi'_{\lv}: g\rightarrow g \cdot  (\o{exp}(\lv)^{-1}) =g \cdot
\o{exp}(-\lv).
\end{equation}

All of  our quantities are conjugation-invariant. Therefore, all
of our calculations reduce to a single hyperbolic element of
$\SLR$,
\begin{equation*}
g= \begin{bmatrix}  e^s & 0 \\ 0 & e^{-s} \end{bmatrix} =
\o{exp}\left( \begin{bmatrix}  s & 0 \\ 0 & -s
\end{bmatrix}\right)
\end{equation*}
whose trace is $\o{Tr}(g) = 2\cosh(s)$. The eigenvalue frame for
the action of $g$ on $\slr$ is
\begin{equation*}
\lx^0_g= \begin{bmatrix}  1 & 0 \\ 0 & -1
\end{bmatrix}, \lx^-_g= \begin{bmatrix}  0 & 0 \\ 1 & 0
\end{bmatrix}, \lx^+_g = \begin{bmatrix}  0 & 1 \\ 0 & 0
\end{bmatrix} ,
\end{equation*}
where
\begin{align*}
g\lx^0_gg^{-1} &= \lx^0_g, \\
g\lx^-_g g^{-1} &= e^{-2s} \lx^-_g  \\
g\lx^+_g g^{-1} &= e^{2s}\lx^+_g .
\end{align*}
Write the vector $\lv\in\slr$ as
\begin{equation*}
 \lv = a \lx^0(g) + b \lx^-(g) + c\lx^+(g)=\begin{bmatrix}  a & b \\ c & -a
\end{bmatrix}.
\end{equation*}

By direct computation, the trace of the induced deformation $\pi_{\lv}(g)$ is
\begin{equation*}
\o{Tr}(\pi_{\lv}(g))= 2\cosh s\cosh \sqrt{a^2+bc} +\frac{2 a \sinh
s \sinh \sqrt{a^2+bc}}{\sqrt{a^2+bc}}
\end{equation*}

Observe that when $\lv =
\begin{bmatrix} 0&b\\c&0 \end{bmatrix}$, which is equivalent to $\ga(\gg)=0$:
\begin{equation*} \o{Tr}(\pi_{\lv}(g))= 2\cosh(s)\cosh(\sqrt{bc})
\end{equation*}
Up to first order, $\o{Tr}(\pi_{\lv}(g))=2\cosh(s)$.

Alternatively,
when $b=c=0$:
\begin{equation*} \o{Tr}(\pi_{\lv}(g))= 2\cosh(s+a)
\end{equation*}
whose Taylor series about $a=0$ does have a linear term. We
assumed that $s>0$, defining our expanding and contracting
eigenvectors. As long as $a>0$, which corresponds to $\ga(\gg)>0$,
the trace of the deformed element $\pi_{\lv}(g)$ is greater than
the original element $g$.

Now consider the deformation $\pi'_{\lv}(g) = g\cdot(\o{exp}(\lv))^{-1}$.   When $b=c=0$:
\begin{equation*} \o{Tr}(\pi'_{\lv}(g))= 2\cosh(s-a)
\end{equation*}
whose Taylor series about $a=0$ has a nonzero linear term. As long
as $a>0$, $\o{Tr}(\pi_{\lv}(g))$ is now less than the original
element $g$. So for this deformation, a positive Margulis
invariant corresponds to a decrease in trace of the original
hyperbolic element.

\begin{lemma}
Consider a hyperbolic $g\in\SLR$, with corresponding closed
geodesic $\partial$ and an affine deformation represented by $\lv
\in \slr$. For the actions of $\lv$ on $\SLR$ by
\begin{itemize}
\item  $\pi_{\lv}(g) = \o{exp}(\lv)\cdot g$ then a positive value
for the Margulis invariant  corresponds to first order lengthening
of $\partial$. \item  $\pi'_{\lv}(g) =  g\cdot\o{exp}(\lv)$ then a
positive value for the Margulis invariant corresponds to first
order shortening of $\partial$.
\end{itemize}
\end{lemma}

\subsection{Deformed parabolic transformations}

As before, we are interested in quantities invariant under
conjugation. Because of this, all of our calculations can be done
with the a very special parabolic transformation in $\SLR$,
\begin{equation*}
p= \begin{bmatrix}  1 & r \\ 0 & 1  \end{bmatrix} = \o{exp} \left(
\begin{bmatrix}  0 & r \\ 0 & 0
\end{bmatrix}  \right)
\end{equation*}
where $r>0$ and whose trace is $\o{Tr}(p) = 2$. We choose a
convenient frame for the action of $p$ on $\slr$:
\begin{equation*}
\lx^u(g)= \begin{bmatrix}  1 & 0 \\ 0 & -1
\end{bmatrix}, \lx^0(g)= \begin{bmatrix}  0 & 0 \\ 1 & 0
\end{bmatrix}, \lx^c(g) = \begin{bmatrix}  0 & 1 \\ 0 & 0
\end{bmatrix} .
\end{equation*}

The trace of the deformation of the element $p$ by the tangent
vector $\lv$ described above is
\begin{equation*}
\o{Tr}(\pi_{\lv}(p))= 2\cosh(\sqrt{a^2+bc})
+\frac{cr}{\sqrt{a^2+bc}}\sinh(\sqrt{a^2+bc})
\end{equation*}

When $\ga(\gg)=0$,
or equivalently when $\lv =
\begin{bmatrix} a&b\\0&-a \end{bmatrix}$, is
\begin{equation*} \o{Tr}(\pi_{\lv}(p))= 2\cosh(a)
\end{equation*}
Thus the trace equals $2$, in terms of $a$, to first order.

Alternatively, when $a=b=0$ in the expression for $\lv$,
\begin{equation*} \o{Tr}(\pi_{\lv}(p))= 2+cr
\end{equation*}
which is linear and increasing in $c$.  As long as $c>0$, which
corresponds to $\ga(\gg)>0$, the trace of the deformed element
$\pi_{\lv}(p)$ majorizes the original element $p$.
\begin{lemma}
Consider a parabolic $g\in\SLR$, and an affine deformation
represented by $\lv \in \slr$. For the actions of $\lv$ on $\SLR$
by
\begin{itemize}
 \item  $\pi_{\lv}(g) = \o{exp}(\lv)\cdot g$ then a positive value for the
 Margulis invariant corresponds to first order increase in the trace of $g$;
\item  $\pi'_{\lv}(g) =  g\cdot\o{exp}(\lv)$ then a positive value for the
 Margulis invariant corresponds to first order decrease in the trace of $g$.
\end{itemize}
\end{lemma}

\end{document}